\documentclass[12pt,a4paper]{article}
\usepackage{amssymb,amsmath,amsthm, amsfonts}
\usepackage{enumerate}
\usepackage[latin1]{inputenc}
\usepackage{bbm}
\usepackage[colorlinks=true,linktocpage=true,linkcolor=blue,citecolor=red]{hyperref}
\usepackage{setspace}

\usepackage[colorlinks=true,linktocpage=true,linkcolor=blue,citecolor=red]{hyperref}
\def\sqr#1#2{{\vcenter{\vbox{\hrule height.#2pt
              \hbox{\vrule width.#2pt height#1pt \kern#1pt \vrule width.#2pt}
              \hrule height.#2pt}}}}

\def \Sum{\displaystyle\sum}

\def \Int{\displaystyle\int}
\def \Frac{\displaystyle\frac}

\def \Sup{\displaystyle\sup}
\def \Lim{\displaystyle\lim}

\def \Max{\displaystyle\max}

\def \ds {\displaystyle}

\def \R{\mathbb{R}}

\def \E{\mathbb{E}}
\def \F{\mathbb{F}}

\def \P{\mathbb{P}}

\def \Ac{{\cal A}}

\def \Fc{{\cal F}}

\def \Ic{{\cal I}}

\def \Lc{{\cal L}}

\def \Sc{{\cal S}}

\def \Xc{{\cal X}}

\def \cd{\cdot}

\def \eps{\varepsilon}

\def \ep{\hbox{ }\hfill$\Box$}

\def\be{\begin{eqnarray}}
\def\ee{\end{eqnarray}}

\def\b*{\begin{eqnarray*}}
\def\e*{\end{eqnarray*}}
\def \nn{\nonumber }
\def \ds{\displaystyle}
\def \1{\mathbbm{1}}

\newtheorem{Lemma}{Lemma}[section]
\newtheorem{Remark}{Remark}[section]

\newtheorem{Theorem}{Theorem}[section]

\newtheorem{Hypothesis}{Hypothesis}[section]

\makeatletter
   
   \@addtoreset{equation}{section}
\makeatother

\def\wt{\widetilde}

\begin{document}

\title{Optimal switching problem and system of  reflected multi-dimensional FBSDEs with random terminal time }
%

\maketitle
\begin{center}

Soufiane Aazizi$^{*}$\footnote{Department of Mathematics, Faculty of
Sciences Semlalia Cadi Ayyad University, B.P. 2390 Marrakesh,
Morocco. Email: {\tt aazizi.soufiane@gmail.com}; {\small\tt
imadefakhouri@gmail.com}\\
$^*$ This author is supported by the Marie Curie Initial Training
Network (ITN) project: ``Deterministic and Stochastic Controlled
Systems and Application", FP7-PEOPLE-2007-1-1-ITN, No. 213841-2.\\
$^{**}$ This author is supported by CNRST ``Centre national pour la recherche scientifique et technique".}
\quad and \quad Imade Fakhouri$^{**1}$ \\ {\it  Université Cadi Ayyad}\\~\\
\end{center}

%

\abstract{In this paper, we study the solvability of a class of multi-dimensional forward backward stochastic differential equations (FBSDEs) with oblique reflection and unbounded stopping time. Under some mild assumptions on the coefficients in such FBSDE, the existence result of adapted solutions is done via a penalization method. The uniqueness is obtained by a verification theorem similarly to the one used by Hu and Tang \cite{HT10}. Finally, we establish the connection with the corresponding optimal switching problem.  This latter is solved by using  the previous results on FBSDEs.}

\medskip

{\bf Key Words. } Backward stochastic differential equations; Oblique reflection; Optimal switching; Unbounded stopping time; Switching problem.

\medskip

\medskip
{\bf AMS Subject Classifications. } 60H10, 93E20

\section{Introduction}
This paper is dedicated to the study of a system of multi-dimensional reflected forward-backward stochastic differential equations (FBSDEs in short) with stopping time not necessarily bounded. In fact, we  generalize the work of Hu and Tang \cite{HT10} to infinite horizon.\\
For $i\in \Lambda:=\{1,\cdots,d\}$ and $t \ge 0$, we define the forward stochastic differential equation (SDE) by
\be \label{SDE}
X_i(t)=x_0+ \int_0^{t\wedge \tau}b(s,X_i(s),i)ds+\int_0^{t\wedge \tau}\sigma(s,X_i(s),i)dW_s,
\ee
and the oblique reflected multi-dimensional backward stochastic differential equation (RBSDE) by
\be\label{RBSDEi}
\left\{\begin{array}{l}
Y_i(t)=g(X_i(\tau))+\Int_{t\wedge \tau}^\tau f(s,X_i(s), Y_i(s),Z_i(s),i) ds+\int_{t\wedge \tau}^\tau dK_i(s)-\int_{t\wedge \tau}^\tau Z_i(s)\, dW(s),\vspace{1mm}\\
Y_i(t)\ge\displaystyle \Max_{j\in \Ic}\{Y_j(t)-C_{i,j}\}, \vspace{1mm}\\
\Int_0^\tau\Big(Y_i(s)-\Max_{j\not=i}\{Y_j(s)-C_{i,j}\}\Big)dK_i(s)=0.
\end{array}
\right.\ee
RBSDEs were firstly studied by El Karoui et al. \cite{E97} for the one dimensional case. Later Gegout-Petit and Pardoux \cite{GP96}  extended this work to the multi-dimensional case with reflection on a boundary convex domain, and recently Hu and Tang \cite{HT10} studied the case of RBSDEs with oblique reflection. In the case of unbounded stopping time, Pardoux \cite{P99} gave existence and uniqueness results of BSDEs under one kind of Lipschitz and monotone assumptions. In the infinite horizon, Hamadène et al. \cite{HLW99}, Akdim and Ouknine \cite{AO06} studied reflected BSDEs and reflected BSDEs with jumps respectively. However, for multi-dimensional reflected  FBSDEs we find only the work of El Asri \cite{E10}, in which the author studied a system of reflected FBSDE  and provided an application to optimal switching problem, but this work suffers from two points: $i)$ The  generator depends only on the forward process. $ii)$ The infinite horizon value of the solution must be zero.

The novelty of this paper lies in the fact that the generator of the BSDE with stopping time depends on the solution $Y_i$ and the process $Z_i$. Here the stopping time is unbounded.  When the stopping time takes infinity, the value of the solution for FBSDE is not necessarily required to be zero. We then prove existence and uniqueness of the solution under one kind of Lipschitz and monotone assumptions. This kind of stopping time will be used to deal with a switching control problem. Given a switching strategy $\alpha \in \Ac$, with $\Ac$ the set of admissible strategies, associated to the controlled process $X^\alpha$ and defined by
\b*
\alpha_t&:=&\sum_{k\geq 0} \zeta_{k}\mathbf{1}_{[\tau_k,\tau_k+1)}(t), \quad t\ge 0,
\e*
here, $\tau_k$ with $k\in \R^+$ is a stopping time such that $\Lim_{k\rightarrow \infty}\tau_k=\tau $ and $\zeta_k$ is an $\Fc_{\tau_k}$-measurable variable with values in $\Lambda$. We consider the total profit at horizon $\tau$ defined by
\b*
J(\alpha_.)=\E^{\alpha_.}\left[ g(X_\tau^{\alpha_.})+\int_0^\tau l(s, X^{\alpha_.}(s), \alpha_s)ds + \Sum_{i\ge 1} C_{\alpha_{i-1}, \alpha_i}  \right],
\e*
where $\E^{\alpha_.}$ is the expectation under probability $P^{\alpha_.}$ defined in (\ref{ProbabilityUnderAlpha}). The optimal switching problem is to maximize the profit $J(\alpha_.)$ with respect to $\alpha_.$, i.e., find an optimal strategy $\alpha_.^*$ such that
\b*
J(\alpha^*_.)=\Sup_{\alpha_. \in \Ac} J(\alpha_.).
\e*

More details on the practical implications of this type of optimal switching problem are given in \cite{CL08} and \cite{PTW09}.

This paper is organized as follows. In Section 2, we state some assumptions and we discuss the case of $X_\tau$ with $\tau$ takes infinity. In Section 3 we prove the existence by a penalization method under one kind of Lipschitz and monotone assumptions, whereas in Section 4 we study the uniqueness via a verification theorem. The last section is devoted to the link between the reflected FBSDEs and the optimal switching problem.

{\bf Notations}. Throughout this paper, we are given a final time $\tau$ which is an $\Fc$-stopping time not necessarily bounded and a probability space $(\Omega, \Fc, P)$
endowed with a $d$ dimensional Brownian motion $W=(W_t)_{t\ge 0}$.  $\{\Fc_t, t\ge 0\}$ is the natural filtration of the Brownian motion augmented by $P$-null sets of $\Fc$.
All the measurability notion will refer to this filtration.
We denote by:\\
{\bf$S^2$} \quad the set of $\R^d$-valued adapted and càdlàg processes $\{Y(t)\}_{t \ge 0}$ such that
$$||Y||_{S^2}:=\E\left[\sup_{0 \le t \le \tau }|Y(t)|^2\right]^{1/2}<+\infty.$$
{\bf $M^2$} \quad  denotes the set of predictable processes $\{Z(t)\}_{t \ge 0}$ with values in $R^{d\times p}$ such that
$$||Z||_{M^2}:=\E\left[\int_0^\tau  |Z(s)|^2ds\right]^{1/2}<+\infty.$$
{\bf $A^2$} \quad is the closed subset of $S^2$ consisting of nondecreasing processes $K=(K_t)_{0 \le t \le \tau}$ with $K_0=0$.\\
{\bf $Q$} \quad the set of process $(y_1, \cdots, y_d)^T \in \R^d$ such that
                \b* y_i &>& y_j-C_{i,j}, \quad \forall i, j \in  \Lambda \mbox{ s.t } i \neq j,  \e*
where $C$ is a real function defined on $\Lambda \times \Lambda$.\\
{\bf $\bar{Q}$} \quad is the closer of domain $Q$ in which the reflected BSDE (\ref{RBSDEi}) evolves, this closer domain is convex and unbounded.

As explained in Hu and Tang \cite{HT10}, each equation of (\ref{RBSDEi}) is independent of others in the interior of $\bar{Q}$ and on bundary $\partial Q$ of domain $Q$ defined by $\partial Q=\ds\cup_{k=1}^d\partial L_k^+$, where for any $k\in \Lambda$
\b*
\partial L_k^+ := \{y\in\R^d: \quad y_k > y_l - C_{k,l}, \mbox{ for any } l\in \Lambda \mbox{ such that } k \neq l\},
\e*
the $k$-th equation is switched to another one, and the solution is reflected along the oblique direction $e_{k}$ which is positive direction of $k$-th coordinate axis.

\section{Preliminaries}
Let us introduce some notations, throughout this paper, we  denote by $\langle \cd , \cd \rangle$ and $|\cd|$ the usual scalar product and the Euclidean
norm for vectors respectively, and by $\|\cd\|$ the trace norm for the matrices. Now, we make the following assumptions:
\vspace{3mm}\\
{\bf(H1)}$\quad$ The functions $b:\R^+\times \Omega \times \R \rightarrow \R$, $\sigma:\R^+\times \Omega \times \R \rightarrow \R$, $g: \R \rightarrow \R$ and
$f:\R^+\times \Omega \times \R \times \R \times \R^p \times \Lambda \rightarrow \R$. Moreover, $b(\cd,x)$, $\sigma(\cd,x)$, $g(x)$
 and $f(\cd,x,y,z,i)$ are all progressively measurable for each $(x,y, z, i) \in \R\times \R \times \R^{p} \times \Lambda$.
 \\
\medskip
{\bf(H2)}$\quad$ $f(\cd, 0, 0, 0):= (f(\cd, 0, 0, 0, 1), \cdots,  f(\cd, 0, 0, 0, d))^T$ belongs to $M^2$.
\\
\medskip
{\bf (H3)}$\quad$ For any $t \ge 0$ , $x, x', y, y', z \in \R$ and $i\in \Lambda$ there exist $\mu_1, \mu_2 \in \R$, $\mu_3 \in \R^+$ and one positive deterministic bounded function $u(t)$, such that
\be
\langle x - x', b(t,x, i)-b(t,x', i) \rangle &\le& \mu_1|x-x'|^2, \quad \P-\mbox{a.s.,}\label{Hbmu1}\\
\langle y - y', f(t,x, y, z, i)-f(t,x, y', z, i) \rangle &\le& \mu_2 u(t)|y-y'|^2, \quad \P-\mbox{a.s.,} \label{Hfmu2}
\ee
and $\int_0^\infty u(t)dt<\infty$, $\int_0^\infty u^2(t)dt<\infty$. \\
{\bf(H4)}$\quad$ For any $t, x, x', y, y', z, z'$  there exist $k\ge 0$  such that
\be
|b(t,x, i)-b(t,x', i)|+|\sigma(t,x, i)-\sigma(t,x', i)| &\le& u(t)|x-x'| \label{Hsigmak1},\\
|f(t,x, y, z, i)-f(t,x', y', z', i)| &\le& u(t)\left(|x-x'|+\|y-y'\|+\|z-z'\| \right). \label{Hfc}
\ee
{\bf(H5)}$\quad$ For any $t, x, x'$  there exist $k_2 \ge 0$ such that
\be
|g(x)-g(x')| &\le& k_2|x-x'| \label{Hgk2}.
\ee
{\bf(H6)}$\quad$ There exist a constant $\lambda \in \R$ such that for any $i \in \Lambda$, a positive constant $C_u$ depending on the function $u$, and $\rho, \eps >0$
\be\begin{array}{l}
\eps^{-1}C_u+2\mu_2u(t)+2\rho^{-1}u^2(t) + 2\eps < \lambda < -2\mu_1 - u(t), \quad t \ge 0,\\
\E\left( \Int_0^\tau e^{\lambda t}(|b(t,0,i)|^2+\|\sigma(t,0,i)\|^2dt\right) < \infty,\label{Hbsigmaffini}\\
\E\left( e^{\lambda \tau}|g(0)|^2+ \Int_0^\tau e^{\lambda t}|f(t,0,0,0,i)|^2dt\right) < \infty.
\end{array}\ee
{\bf(H7)}$\quad$ For any $i\in \Lambda$ and $\tau \in [0, +\infty]$ we have
\be
\E\left(\int_0^\tau|b(s,0,i)|ds\right)^2 + \E\int_0^\tau|\sigma(s,0,i)|^2ds &<& \infty.\label{H b+sigma}
\ee
\begin{Remark}
For simplicity, we take the same function $u(t)$ in (\ref{Hfmu2}), (\ref{Hfc}) and (\ref{Hbsigmaffini}).
\end{Remark}
The reflected BSDE (\ref{RBSDEi}) evolves in the closure $\bar{Q}$ of domain $Q$.
As a preparation, we first recall a lemma which is proved by Yin \cite{Y08}:
\begin{Lemma} \label{Lemma SDE}(See \cite[Remark 2.1 and Lemma 3.2]{Y08})\\
Assume (\ref{Hbmu1}), (\ref{Hsigmak1}) and (\ref{Hbsigmaffini}) hold, where $\lambda< -2\mu_1-k_1$. Then the forward SDE (\ref{SDE}) admits a unique solution $\{X(t)\}_{t \ge 0}$ satisfying
\b*
\E\left( \sup_{0 \le t \le \tau}e^{\lambda t}|X(t)|^2+\int_0^\tau e^{\lambda t}|X(t)|dt\right) &<&\infty.
\e*

\end{Lemma}

Before proving existence, we shall discuss the case of $X_\tau$ with $\{\tau=+\infty\}$ which appears in the BSDE (\ref{SDE}).
Under Hypothesis (\ref{Hsigmak1}) and (\ref{H b+sigma}), the integral $\Int_0^\infty \sigma(s,X_i(s),i)dW_s$ is well defined and is an $L^2$-bounded  martingale. Thus, it is easy to show that
\b*
\Lim_{t\rightarrow \infty} \E\left[ \left|\int_0^{t} \sigma(s,X_i(s),i)dW_s-\int_0^\infty \sigma(s,X_i(s),i)dW_s\right|^2\right]=0.
\e*
Now, we define
\b*
\mathcal{X}=x_0+ \int_0^\infty b(s,X_i(s),i)ds+\int_0^\infty \sigma(s,X_i(s),i)dW_s.
\e*
Then from (\ref{SDE})  we have
\b*
\Xc - X_\tau= \int_\tau^\infty b(s,X_i(s),i)ds+\int_\tau^\infty \sigma(s,X_i(s),i)dW_s, \quad \forall t\ge 0.
\e*
It is obvious that $ \lim_{\tau \rightarrow \infty}\E|\Xc - X_\tau|^2 =0$, so that $\mathcal{X}=\lim_{\tau \rightarrow \infty}X_\tau$ in $L^2$ and we denote it by $X_\infty$.\\
For more details on the process $X_\tau$ with $\tau \in [0, \infty]$, we send the reader to \cite{W04}.

\section{Existence}

In this section, we shall prove an existence theorem of solution of FBSDE (\ref{SDE})-(\ref{RBSDEi}). Our setup contains the case $\{\tau \equiv +\infty\}$ as a particular case. Let us firstly make the following assumptions on the cost function $C$ which are standard in the optimal switching problem.

\begin{Hypothesis}\label{k} $ $
\begin{enumerate}[(i)]
  \item For any $(i,j)\in \Lambda\times \Lambda$, $C_{i,j}\ge 0$.

  \item For any $(i,j,l)\in \Lambda\times \Lambda\times \Lambda$, such that $i\neq j$ and $j \neq l$, we have
\b*C_{i,j}+C_{j,l} &\ge& C_{i,l}.\e*
\end{enumerate}
\end{Hypothesis}

For $n \ge 0$, let us introduce  the following penalized BSDE for any $t \ge 0$ and $i\in \Lambda$:
\be
\label{penalized}\begin{array}{rcl}
Y_i^n(t)&=&g(X_i(\tau)) +\Int_{t\wedge \tau}^\tau f(s,X_i(s),Y_i^n(s),Z_i^n(s),i)-\Int_{t\wedge \tau}^\tau Z_i^n(s)\ dW(s)\\
&&+n\Sum_{l=1}^d\Int_{t\wedge \tau}^\tau\big(Y_i^n(s)-Y_l^n(s)+C_{i,l}\big)^-\,ds.
\end{array}\label{PBSDE}
\ee
Define
\b*
K^n_t=n\sum_{l=1}^d\int_0^{t\wedge \tau}\big(Y_i^n(s)-Y_l^n(s)+C_{i,l}\big)^-\,ds
\e*
Note that when $l=i$, we have
\be(Y_i^n(s)-Y_l^n(s)+C_{i,l})^-=0. \ee

From the classical result of Chen \cite{C98}, for any
$n \ge 0$, BSDE (\ref{penalized}) has a unique solution $(Y^n,Z^n)$ in
the space $S^2\times M^2$.

We're going to prove that the triplet $(Y^n, Z^n, K^n)$ converges to the solution of RBSDE (\ref{RBSDEi}). To do so, we first need the following a priori estimation.

\subsection{A priori estimation }
In this subsection, we derive two lemmas on the a priori estimation of the penalized BSDE (\ref{PBSDE}), which will play a primordial role in the sequence $(Y^n, Z^n, K^n)$ convergence proof.
\begin{Lemma}\label{aprioriK} Let the Hypotheses ($H2$), (\ref{Hfmu2}), (\ref{Hfc}) and (\ref{Hbsigmaffini}) hold true
 and assume that  $\forall i \in \Lambda$, $g(X_i(\tau))\in L^2(\Omega, {\cal F}_\tau, P, \R^d)$ takes values in
$\bar{Q}$. Then there exists a constant $C>0$ (independent of $n$), such that
\be\begin{array}{rcl}
&&\E\left(\Sup_{0\le t \le \tau }e^{\lambda t}\Big|(Y_i^n(t)-Y_j^n(t)+C_{i,j})^-\Big|^2+ n^2\Int_0^\tau e^{\lambda t}\Big|(Y_i^n(t)-Y_j^n(t)+C_{i,j})^-\Big|^2\, dt \right)\\
&\le& \\
&&C_u \E\left(\Int_0^\tau e^{\lambda s}(|f(s,0,0,0)|^2+|X_i(s)|^2+|X_i(s)-X_j(s)|^2)+|Y_i^n(s)|^2+|Z_i^n(s)|^2 \, ds\right),
\end{array}
\ee
where $C_u$ is a constant depending on $u(t)$.
\end{Lemma}
{\bf Proof. }\vspace{2mm} \\
For simplicity, denote $\bar{Y}_{ij}^n= Y_i^n(t)-Y_j^n(t)+C_{i,j}$, we have for all $t \ge 0$, $ i\in \Lambda$
\be\begin{array}{rcl}
\bar{Y}_{ij}^n(t)&=&\bar{Y}_{ij}^n(T)+\Int_{t\wedge \tau}^\tau \left[f(s,X_i(s),Y_i^n(s),Z_i^n(s),i)-f(s,X_j(s),Y_j^n(s),Z_j^n(s),j)\right]\, ds\\
[0.5cm]&&+n\Sum_{l=1}^d\int_{t\wedge \tau}^\tau\bar{Y}_{il}^n(s)^-\,ds-n\sum_{l=1}^d\int_{t\wedge \tau}^\tau\bar{Y}_{jl}^n(s)^- -\Int_{t\wedge \tau}^\tau \left[Z_i^n(s) -Z_j^n(s) \right]\ dW(s).
\end{array}
\ee
For $i,j\in\Lambda$, if we denote $L_{ij}^n$ the local time  of the semi-martingale  $\bar{Y}_{ij}^n(t)$, then we get by Tanaka  formula
\b*\begin{array}{rcl}
&&\bar{Y}_{ij}^n(t)^- +n\Sum_{l=1}^d\int_{t\wedge \tau}^\tau I_{\Lc_{ij,n}}(s)\left[\bar{Y}_{il}^n(s)^- - \bar{Y}_{jl}^n(s)^-\right]\,ds+{\frac{1}{2}}\int_{t\wedge \tau}^\tau \, dL_{ij}^n(s)\\
 &&=\Int_{t\wedge \tau}^\tau  I_{\Lc_{ij,n}}(s)\left[f(s,X_i(s),Y_i^n(s),Z_i^n(s),i)-f(s,X_j(s),Y_j^n(s),Z_j^n(s),j)\right]\, ds\\
&&-\Int_{t\wedge \tau}^\tau  I_{\Lc_{ij,n}}(s) \left[Z_i^n(s) -Z_j^n(s) \right]\ dW(s),
\end{array}
\e*
where for $i,j\in\Lambda$,
\be\begin{array}{c}
\Lc_{ij,n}:=\{(s,\omega):\bar{Y}_{ij}^n(s)<0\}. \end{array}
\ee
Applying Itô's formula for $e^{\lambda t\wedge \tau} |\bar{Y}_{ij}^n(t)^-|^2$ yields
\be\begin{array}{rcl}
&&e^{\lambda t} |\bar{Y}_{ij}^n(t)^-|^2 +(2n+\lambda)\Int_{t\wedge \tau}^\tau e^{\lambda s}|\bar{Y}_{ij}^n(s)^-|^2 \, ds+\Int_{t\wedge \tau}^\tau I_{\Lc_{ij,n}}(s)e^{\lambda s}|Z_i^n(s) -Z_j^n(s)|^2 \, ds\\
&=&2\Int_{t\wedge \tau}^\tau I_{\Lc_{ij,n}}(s)e^{\lambda s}\bar{Y}_{ij}^n(s)^- \left[f(s,X_i(s),Y_i^n(s),Z_i^n(s),i)-f(s,X_j(s),Y_j^n(s),Z_j^n(s),j)\right]\, ds\\
&&-2\Int_{t\wedge \tau}^\tau I_{\Lc_{ij,n}}(s) e^{\lambda s}\bar{Y}_{ij}^n(s)^-\left[Z_i^n(s) -Z_j^n(s) \right]\ dW(s)\\
&&+2n\Int_{t\wedge \tau}^\tau e^{\lambda s}\bar{Y}_{ij}^n(s)^-\bar{Y}_{ji}^n(s)^- \,ds+2n\Sum_{l\neq i, l\neq j}\Int_{t\wedge \tau}^\tau e^{\lambda s}\bar{Y}_{ij}^n(s)^- \left[\bar{Y}_{jl}^n(s)^- -\bar{Y}_{il}^n(s)^- \right]\,ds,\label{Tanaka}
\end{array}
\ee
since we have
\b* \int_{t\wedge \tau}^\tau \bar{Y}_{ij}^n(s)^-\,dL_{ij}^n(s)=0,\quad \forall t\ge0.
\e*
From other side, since $C_{i,j}+C_{j,i} \ge 0,$ then  $ \bar{Y}_{ji}^n(s)^-\bar{Y}_{ij}^n(s)^- =0$. In fact
\b*
\{y\in \R^d: y-y'+C_{i,j}<0\} \cap \{y\in \R^d: y'-y+C_{j,i}<0\}&=&\emptyset.
\e*
Also we know that for two real numbers $x_1$ and $x_2$, we have $x_1^--x_2^-\le (x_1-x_2)^-$, then
\b* \begin{array}{rcl}
I_{\Lc_{ij,n}}(s)\left[\bar{Y}_{jl}^n(s)^- -\bar{Y}_{il}^n(s)^-\right] &\le& I_{\Lc_{ij,n}}(s)\left[\bar{Y}_{jl}^n(s) -\bar{Y}_{il}^n(s)\right]^-\\
&=& I_{\Lc_{ij,n}}(s)(Y_j^n(s)-Y_i^n(s)+C_{j,l}-C_{i,l})^-=0,
\end{array}
\e*
in view that
\b*
\{y \in \R^d:y-y'+C_{i,j}<0\}\cap\{y \in \R^d: y'-y+C_{j,l}-C_{i,l}<0\}=\emptyset.
\e*
Combining this together with (\ref{Tanaka}) and taking expectation, we get
\be\label{TanakaE}\begin{array}{rcl}
&&\E \left(e^{\lambda t \wedge \tau}|\bar{Y}_{ij}^n(t)^-|^2\right) +(2n+\lambda)\E\Int_{t\wedge \tau}^\tau e^{\lambda s}|\bar{Y}_{ij}^n(s)^-|^2 \, ds\\
&&\quad +\E\Int_{t\wedge \tau}^\tau I_{\Lc_{ij,n}}(s)e^{\lambda s}|Z_i^n(s) -Z_j^n(s)|^2 \, ds\\
&&\le 2\E\Int_{t\wedge \tau}^\tau e^{\lambda s} \bar{Y}_{ij}^n(s)^-|f(s,X_i(s),Y_i^n(s), Z_i^n(s),i)- f(s,X_j(s),Y_j^n(s), Z_j^n(s),j)|\,ds\\
%
&&\le 2\E\Int_{t\wedge \tau}^\tau e^{\lambda s} \bar{Y}_{ij}^n(s)^-\Big[|f(s,X_i(s),Y_i^n(s), Z_i^n(s),i)- f(s,X_i(s),Y_i^n(s), Z_i^n(s),j)|\\
&& \hspace{30mm}+\left.|f(s,X_i(s),Y_i^n(s), Z_i^n(s),j)- f(s,X_j(s),Y_j^n(s), Z_j^n(s),j)|\right]\,ds\\
&&\le 2\E\Int_{t\wedge \tau}^\tau  e^{\lambda s}u(s)\bar{Y}_{ij}^n(s)^-\Big[u^{-1}(s)f(s,0,0,0)+|X_i(s)|+|Y_i^n(s)|+|Z_i^n(s)|\\
&& \hspace{35mm}+|X_i(s)-X_j(s)|+|\bar{Y}_{ij}^n(s)|+|Z_i^n(s)-Z_j^n(s)|\Big]ds\\
&&\le \E\Int_{t\wedge \tau}^\tau e^{\lambda s}(1+u(s)+5u^2(s))|\bar{Y}_{ij}^n(s)^-|^2 \, ds \\
&&+ \Frac12 \,\E\Int_{t\wedge \tau}^\tau I_{\Lc_{ij,n}}(s)e^{\lambda s} \Big[|f(s,0,0,0)|^2+|X(s)|^2+|Y_i^n(s)|^2+|Z_i^n(s)|^2 \\
&& \hspace{35mm}+|X_i(s)-X_j(s)|+|Z_i^n(s)-Z_j^n(s)|^2 \,\Big]   ds.\\
\end{array} \ee

Applying Gronwall's inequality and  Lemma \ref{Lemma SDE}, it follows that
\b*&&\E \left(e^{\lambda t\wedge \tau}|\bar{Y}_{ij}^n(t)^-|^2\right)\\
&\le& C_u\E\Int_0^\tau I_{\Lc_{ij,n}}(s)e^{\lambda s}\Big[|f(s,0,0,0)|^2+|X_i(s)|^2+|X_i(s)-X_j(s)|^2+|Y_i^n(s)|^2+|Z_i^n(s)|^2 \Big]\, ds ,
\e*
and
\b*
&&n\E\Int_{t\wedge \tau}^\tau I_{\Lc_{ij,n}}(s)e^{\lambda s}|\bar{Y}_{ij}^n(s)^-|^2 \, ds+\E\Int_{t\wedge \tau}^\tau I_{\Lc_{ij,n}}(s)e^{\lambda s}|Z_i^n(s) -Z_j^n(s)|^2 \, ds\\
&&\le  C_u\E\Int_0^\tau I_{\Lc_{ij,n}}(s)e^{\lambda s}\Big[|f(s,0,0,0)|^2+|X_i(s)|^2+|X_i(s)-X_j(s)|^2+|Y_i^n(s)|^2+|Z_i^n(s)|^2 \Big]\, ds ,\e*
where $C_u$ is a constant depending on $u(t)$ and that will play a crucial role in Lemma \ref{apriori} below. It then follows from Burkholder-Davis-Gundy's inequality applied  to  (\ref{Tanaka})
\b*
\displaystyle &&\E\left[\sup_{0\le t \le \tau}e^{\lambda t}|\bar{Y}_{ij}^n(t)^-|^2\right]\\
&&\le C_u\E\Int_0^\tau I_{\Lc_{ij,n}}(s)e^{\lambda s}\Big[|f(s,0,0,0)|^2+|X_i(s)|^2+|X_i(s)-X_j(s)|^2 +|Y_i^n(s)|^2+|Z_i^n(s)|^2 \Big]\, ds.\nn
\e*
Now from the first inequality in  (\ref{TanakaE}), we get
\b*
&& (2n+\lambda)\E\Int_0^\tau I_{\Lc_{ij,n}}(s)e^{\lambda s}|\bar{Y}_{ij}^n(s)^-|^2 \, ds\\
&&\le  \left(n+\frac{C_u}{n}\right)\E\int_0^\tau e^{\lambda s}|\bar{Y}_{ij}^n(s)^-|^2\, ds\\
&&\,+{\frac{C_u}{n}}\E\int_0^\tau  I_{\Lc_{ij,n}}(s)e^{\lambda s}\Big[|f(s,0,0,0)|^2+|X_i(s)|^2+|X_i(s)-X_j(s)|^2 +|Y_i^n(s)|^2+|Z_i^n(s)|^2 \Big]\, ds .
\e*
For  $n$ large enough we finally deduce that
 \b*
 &&\displaystyle n^2\E\int_0^\tau e^{\lambda s} |\bar{Y}_{ij}^n(s)^-|^2\, ds\nn \\
&\le& C_u\E\Int_0^\tau I_{\Lc_{ij,n}}(s)e^{\lambda s}\Big[|f(s,0,0,0)|^2+|X_i(s)|^2+|X_i(s)-X_j(s)|^2 +|Y_i^n(s)|^2+|Z_i^n(s)|^2 \Big]\, ds.
\e*
\ep

Then, we are able to prove the following estimation:
\begin{Lemma}\label{apriori}Assume ($H2$), (\ref{Hfmu2}), (\ref{Hfc}), (\ref{Hgk2}) and (\ref{Hbsigmaffini}) hold true. Let us also assume that $\forall i \in \Lambda$, $g(X_i(\tau))\in L^2(\Omega,{\cal F}_\tau,P;\R^d)$ takes values in $\bar{Q}$. Then there exists a constant $C>0$, such that
\be\begin{array}{rcl}\label{estimateyz}
&&\E\left(\Sup_{0\le t \le \tau }e^{\lambda t}|Y_i^n(t)^-|^2+\Int_0^\tau e^{\lambda t}|Y_i^n(t)|^2dt+ \Int_0^\tau e^{\lambda t}|Z_i^n(t)|^2dt\right)\\
&&\le\\
&& C \E\left(e^{\lambda \tau}|g(0)|^2+\Sup_{0\le t\le\tau}e^{\lambda t}|X_i(t)|^2 + \Int_0^\tau e^{\lambda s}\Big[|f(s,0,0,0)|^2+|X_i(s)|^2+ \Sum_{j=1}^d|X_i(s)-X_j(s)|^2  \Big]\, ds\right),
\end{array}
\ee
where $\lambda > \eps^{-1}C_u+2\mu_2u(t)+2\rho^{-1}u^2(t) + 2\eps $ and C depends only on $k_2$, $\eps$, $\rho$ and the function $u$.
\end{Lemma}
\newpage
{\bf Proof. }\vspace{2mm} \\
Applying Itô's formula to $e^{\lambda t\wedge \tau}|Y_i^n(t)|^2$, we obtain:
\begin{eqnarray}\label{finalito}
& &e^{\lambda t}|Y_i^n(t)|^2+\int_{t\wedge \tau}^\tau e^{\lambda s}(\lambda |Y_i^n(s)|^2+|Z_i^n(s)|^2)\, ds\nonumber \\
&=& e^{\lambda \tau}|g(X_i(\tau))|^2+2\int_{t\wedge \tau}^\tau e^{\lambda s}Y_i^n(s)\cdot\Big[ f(s,X_i(s),Y_i^n(s),Z_i^n(s),i)+n\sum_{l=1}^d (Y_i^n(s)-Y_l^n(s)+C_{i,l})^-\Big]\, ds\nonumber\\
 & &-2\int_{t\wedge \tau}^\tau e^{\lambda s}Z_i(s)dW(s).
\end{eqnarray}
Taking expectation and using the fact that for any arbitrary $\eps>0$ and any $\rho<1$ arbitrarily close to one,
\b*
2\langle y, f(t,x,y,z)\rangle \le (2\mu_2u(t)+2\rho^{-1}u^2(t) + \eps)|y|^2 + \rho\|x\|^2+\rho\|z\|^2+\eps^{-1}|f(t,
,0,0)|^2,\e*
combined with Lemma \ref{aprioriK} we get
\b*
 & &e^{\lambda t\wedge \tau}\E|Y_i^n(t)|^2+\E\int_{t\wedge \tau}^\tau e^{\lambda s}(\lambda |Y_i^n(s)|^2+\rho|Z_i^n(s)|^2)\, ds\\
&\le&\E\left[e^{\lambda \tau}|g(0)|^2\right]+ k_2 \E\sup_{0\le t\le\tau}e^{\lambda t}|X_i(t)|^2+ \E\int_{t\wedge \tau}^\tau (2\mu_2u(t)+2\rho^{-1}u^2(t) + 2\eps)\,e^{\lambda s}|Y_i^n(s)|^2\, ds  \\
&& +\eps^{-1} \E\int_{t\wedge \tau}^\tau e^{\lambda s} | f(s,0,0,0)|^2\, ds +\rho \E\int_{t\wedge \tau}^\tau e^{\lambda s}\Big[|X_i(s)|^2+|Z_i(s)|^2\Big]\,ds\\
& &+n^2\eps^{-1}\E\int_{t\wedge \tau}^\tau e^{\lambda s}\sum_{l=1}^d\left((Y_i^n(s)-Y_l^n(s)+C_{i,l})^-\right)^2\, ds.
\e*
For $\bar{\lambda}:= \lambda-\eps^{-1}C_u-2\mu_2u(t)-2\rho^{-1}u^2(t) - 2\eps > 0$ and $\bar{\rho}=1-\rho-\eps^{-1}C_u >0$, we get
\b*
 & &e^{\lambda t}\E|Y_i^n(t)|^2+\E\int_{t\wedge \tau}^\tau e^{\lambda s}(\bar{\lambda} |Y_i^n(s)|^2+\bar{\rho}|Z_i^n(s)|^2)\, ds\\
&\le&C \E\left(e^{\lambda \tau}|g(0)|^2+\sup_{0\le t\le\tau}e^{\lambda t}|X_i(t)|^2 + \Int_0^\tau e^{\lambda s}\Big[|f(s,0,0,0)|^2+|X_i(s)|^2+ \Sum_{j=1}^d|X_i(s)-X_j(s)|^2  \Big]\, ds\right),
\e*
where $C$ depends on $k_2$, $\eps$, $\rho$ and the function $u$.
Finally, we deduce by an argument already used. This completes the proof.
\ep
\medskip

This will allow us to prove the convergence of the sequence $(Y^n, Z^n, K^n)$.
\medskip
\subsection{Convergence of the sequence $(Y^n, Z^n, K^n)$}
Now, we will prove that $(Y^n, Z^n, K^n)_{n\ge 0}$ is a Cauchy sequence.
\begin{Lemma} \label{LemmaCauchy}
The sequence $\{(Y^n,Z^n)\}_n$ is a Cauchy sequence in the space $S^2\times M^2$.
\end{Lemma}
{\bf Proof. }\vspace{2mm} Denote: \\
$Y_i^n(t)-Y_i^m(t)=\bar{Y}_i^{n,m}(t)$, \\
$Z_i^n(t)-Z_i^m(t)=\bar{Z}_i^{n,m}(t)$, \\
$\bar{Y}_{ij}^n= Y_i^n(t)-Y_j^n(t)+C_{i,j}$. \\

Applying Itô's formula to $e^{\lambda t\wedge \tau}|\bar{Y}_i^{n,m}(t)|^2$, we have for $i\in \Lambda$

\begin{eqnarray}\label{itoynm}
& &\E\left(e^{\lambda t\wedge \tau}|\bar{Y}_i^{n,m}(t)|^2 \right)+ \E\int_{t\wedge \tau}^\tau e^{\lambda s}(\lambda |\bar{Y}_i^{n,m}(s)|^2+|\bar{Z}_i^{n,m}(s))|^2\, ds\nonumber\\
&=&2\E\int_0^\tau e^{\lambda s}\bar{Y}_i^{n,m}(s)(f(s,X_i(s),Y_i^n(s), Z_i^n(s),i)- f(s,X_i(s),Y_i^m(s),Z_i^m(s),i))\, ds\nonumber \\
& & +2n\E\int_{t\wedge \tau}^\tau e^{\lambda s}\bar{Y}_i^{n,m}(s)\bar{Y}_{ij}^n(s)^-\, ds -2m\E\int_{t\wedge \tau}^\tau e^{\lambda s}\bar{Y}_i^{n,m}(s)\bar{Y}_{ij}^m(s)^-\, ds\nonumber\\
&\le &C_\alpha\E\int_0^\tau u^2(s)e^{\lambda s}|\bar{Y}_i^{n,m}(s)|^2\, ds+\alpha\E\int_0^\tau e^{\lambda s}|\bar{Z}_i^{n,m}(s)|^2\, ds\nonumber \\
& & +2\E\int_{t\wedge \tau}^\tau n^2e^{\lambda s}|\bar{Y}_{ij}^n(s)^-|^2\, ds +m^2e^{\lambda s}|\bar{Y}_{ij}^m(s)^-|^2\,ds.
\end{eqnarray}
From Lemma \ref{aprioriK},  and by applying Gronwall's Lemma for  $\alpha <1$. We obtain
\b*
\forall m\ge n, \quad \quad \sup_{0 \le t \le \tau}\E\left(e^{\lambda t}|\bar{Y}_i^{n,m}(t)|^2 \right) &\le& \frac{C}{n}.
\e*
We deduce also
\b*
\forall m\ge n, \quad \quad \E \int_0^\tau e^{\lambda t}|\bar{Z}_i^{n,m}(t)|^2 dt &\le& \frac{C}{n}.
\e*

%
We rewrite again Itô's formula for  $e^{\lambda t}|\bar{Y}_i^{n,m}(t)|^2$, using then Burkholder-Davis-Gundy's inequality and some argument already used, we obtain for
$i\in\Lambda$,
\b*
\E\sup_{0\le t\le \tau}e^{\lambda t}|\bar{Y}_i^{n,m}(t)|^2 &\le& \Frac{C}{n}.
\e*
\ep \\

Let us now define the process $Y_t = \lim_{n\rightarrow +\infty} Y^n_t$ in the sense that $Y^n$ converges to $Y$ in $S^2$, and $Z_t = \lim_{n\rightarrow +\infty} Z^n_t$ in the sense that $Z^n$ converges to $Z$ in $M^2$.\\
We define also:
\be\label{Kn}
K_i^n(t):=n\int_0^{t\wedge \tau}\sum_{l=1}^d(Y_i^n(s)-Y_l^n(s)+C_{i,l})^-\, ds, \quad i\in \Lambda.
\ee
From the expression of BSDE (\ref{penalized}), we have
\begin{equation}\label{Kin}
K_i^n(t)=Y_i^n(t)-Y_i^n(0)+K_i^n(\tau)+\int_0^{t\wedge \tau}  f(s,X_i(s),Y_i^n(s),Z^n_i(s),i)ds-\int_0^{t\wedge \tau} Z^n_i(s)dW(s),\quad i\in\Lambda.
\end{equation}
Set
\begin{equation}\label{eqK}
K_i(t):=Y_i(t)-Y_i(0)+\int_0^{t\wedge \tau}  f(s,X_i(s),Y_i(s),Z_i(s),i)ds-\int_0^{t\wedge \tau} Z_i(s)dW(s),\quad i\in\Lambda.
\end{equation}
Then, we deduce immediately that $K^n$ converges to $K$ in $S^2$.\\
 Finally, it remains to show that
 \be \label{limitK}
 \int_0^\tau\left(Y_i(s)-\Max_{j\not=i}[Y_j(s)-C_{i,j}]\right)^+\, dK^n_i(s)=0,  \quad i\in \Lambda.
 \ee
However, we have from (\ref{Kn}) that for $i\in \Lambda$

 \b*
&& \int_0^\tau\left(Y^n_i(s)-\Max_{j\not=i}[Y^n_j(s)-C_{i,j}]\right)^+ dK_i^n(s)\\
&=& n\sum_{l=1}^d\int_0^\tau\left(Y^n_i(s)-\Max_{j\not=i}[Y^n_j(s)-C_{i,j}]\right)^+ (Y_i^n(s)-Y_l^n(s)+C_{i,l})^-\, ds,
 \e*
which is equal to zero by construction, then as $n \rightarrow \infty$, from \cite[Lemma 5.8]{GP96} we have $(\ref{limitK}$).
In fact, we have shown the existence of the solution of the reflected BSDEs (\ref{RBSDEi}):

\begin{Theorem}\label{existence}
Let the Hypotheses ($H1-H7$) hold. Assume that $g(X(\tau)) \in L^2(\Omega,{\cal F}_\tau,P;\R^d)$ takes values in $\bar{Q}$. Then RBSDE (\ref{RBSDEi}) has a solution $(Y,Z,K)$ in $S^2\times
M^2\times A^2$.
\end{Theorem}
%
%

\medskip

\section{Verification theorem}

A switching strategy $\alpha$ consist in a sequence $\alpha:= (\tau_k, \zeta_k)_{k\ge 1}$, where $(\tau_k)_{k\ge 1}$ is an increasing sequence of $\F$-stopping times smaller than $\tau$, and $\zeta_k$ are $\Fc_{\tau_k}$-measurable random variables valued in $\Lambda$. For an initial regime $i_0$ we define an admissible strategy as follows:
\be\label{admi}
\alpha_t &:=& \Sum_{k\ge 0} \zeta_k1_{[\tau_k, \tau_{k+1} ]}(t), \quad \quad t \ge 0,
\ee
with $\tau_0=0$ and $\zeta_0=i_0$.\\
We denote by $\Ac(t)$ the set of admissible strategies  starting at time $t$ and $\Ac_i(t)$ the subset of $\Ac(t)$ starting at time $t$ from the mode $i$
\b*
\Ac_i(t)&:=& \{\alpha \in \Ac(t): \quad \alpha_t =i\}.
\e*
For any $\alpha_\cdot$ we define the process $A^{\alpha_\cdot}$ by
\be\label{UA}
A^{\alpha_\cdot}(s)=\Sum_{k\ge0} C_{\zeta_k, \zeta_{k+1}} 1_{[\tau_k, \tau]}(s).
\ee
Given a strategy $\alpha \in \Ac$ we define the following BSDE:
\be\label{UV}
U(s) =g(X^{\alpha_\cd}(\tau)) +A^{\alpha_\cd}(\tau)-A^{\alpha_\cd}(s)+\int_{s\wedge \tau}^\tau\psi(r,U(r),V(r),\alpha(r)) dr -\int_{s\wedge \tau}^\tau V(r)dW(r),\quad s \ge t.\ee
\medskip
This BSDE has a solution in $S^2 \times M^2$ denoted $(U^{\alpha_\cd}, V^{\alpha_\cd})$, to prove this, it is enough to write for $s \ge t$
\b*
\tilde{U}(s)&=&U(s) +  A^{\alpha_\cd}(s),\\
\tilde{V}(s)&=&V(s).
\e*
Then we get  from (\ref{UV})
\b*
\tilde{U}(s) &=&g(X^{\alpha_\cd}(\tau)) +A^{\alpha_\cd}(\tau)+\int_{s\wedge \tau}^\tau\psi(r,\tilde{U}(r)-  A^{\alpha_\cd}(r),\tilde{V}(r),\alpha(r)) dr -\int_{s\wedge \tau}^\tau V(r)dW(r),\quad \quad s \ge t.
\e*
Which has solution from standard arguments.
We impose the following stronger assumptions:
\begin{Hypothesis}\label{k'}
(i) For any $(i,j)\in \Lambda\times \Lambda$, $C_{i,j}\ge 0$.

(ii) For any $(i,j,l)\in \Lambda\times\Lambda\times \Lambda$ such
that $i\not=j$ and $j\not=l$,
$$C_{i,j}+C_{j,l}> C_{i,l}.$$
\end{Hypothesis}

With the following representation of the solution of BSDE  (\ref{RBSDEi}), we have immediately
the uniqueness of the solution.

\begin{Theorem}\label{uniqueness} Let us suppose that the Hypotheses (H2), (H4) and \ref{k'} hold. Let us also assume that
$g(X(\tau)) \in L^2(\Omega,{\cal F}_\tau,P;\R^d)$ takes values in $\bar{Q}$.
Let $(\wt Y,\wt Z,\wt K)$ be a solution in $(S^2,M^2,K^2)$ to
RBSDE (\ref{RBSDEi}). Then
\begin{enumerate}[(i)]
  \item For any $\alpha(\cdot)\in {\cal A}_i(t)$, we have:
\begin{equation}\label{inequality}\wt Y_i(t)\le U^{\alpha(\cdot)}(t), \quad P-a.s.
\end{equation}

  \item Set $\tau^*_0=t$, $\zeta^*_0=i$ and define the sequence
$\{\tau_j^*,\zeta_j^*\}_{j=1}^{\infty}$ in an inductive way as
follows:
\begin{equation}\label{theta*}
\tau^*_j:=\inf\{s\ge \tau^*_{j-1}: \tilde{Y}_{\zeta^*_{j-1}}(s)=\max_{l\not={\zeta^*_{j-1}}}\{\tilde{Y}_l(s)-C_{\zeta^*_{j-1},l}\}\wedge \tau,
\end{equation}
and $\zeta_j^*$ is ${\cal F}_{\tau_j^*}$-measurable random
variable such that
$$\wt Y_{\zeta^*_{j-1}}(\tau_j^*)=\wt Y_{\zeta_j^*}(\tau^*_j)-C_{\zeta_{j-1}^*,\zeta^*_j},$$
with $j=1,2,\cdots.$

Then, the following  switching strategy:
\begin{equation}
\alpha^*_s=i\1_{\{t\}}(s)+\sum_{j\ge1} \zeta_{j-1}^*\1_{(\tau^*_{j-1},\tau^*_j]}(s),
\end{equation}
is admissible, i.e., $\alpha^*_\cdot \in {\cal A}_i(t)$ and we have,
$$\wt Y_i(t)=U^{a^*(\cdot)}(t).$$
Moreover, $\wt Y(t)$:
$$\wt Y_i(t)=\mathop{ess \sup}_{\alpha_\cdot \in {\cal A}_t^i} U^{\alpha_\cdot}(t), \quad i\in \Lambda, \quad t \ge 0.$$
\end{enumerate}
Therefore RBSDE (\ref{RBSDEi}) has a unique solution.
\end{Theorem}

{\bf Proof.}   We prove w.l.o.g (i) and (ii) for the particular  case of $t=0$.

(i) We define
\be
\wt Y^{\alpha_\cd}(s)&=&\mathop{\sum_{i\ge1}}_{{\tau_i\in[0, \tau)}}\wt Y_{\zeta_{i-1}}(s)\1_{[\tau_{i-1}, \tau_i)}(s)+g(X^{\alpha_\tau}(\tau))\1_{\{\tau\}}(s),\label{UY}\\
\wt Z^{\alpha_\cd}(s)&=&\mathop{\sum_{i\ge1}}_{{\tau_i\in[0, \tau)}}\wt Z_{\zeta_{i-1}}(s)\1_{[\tau_{i-1}, \tau_i)}(s),\label{UZ}\\
\wt K^{\alpha(\cdot)}(s)&=&\mathop{\sum_{i\ge1}}_{{\tau_i\in[0, \tau)}}\int_{\tau_{i-1}\wedge s}^{\tau_i\wedge s}d\wt K_{\zeta_{i-1}}(r).\label{UK}
 \end{eqnarray}

The process  $\wt Y^{\alpha_\cd}(\cd)$ is càdlàg with jump $\wt Y_{\alpha_i}(\tau_i)-\wt Y_{\alpha_{i-1}}(\tau_i)$ at $\tau_i \in [0, \tau]$, $i \in \Lambda$,  it follows that

\begin{eqnarray*}
\wt Y^{\alpha_\cd}(s)-\wt Y^{\alpha_\cd}(0)&=&\mathop{\sum_{i\ge1}}_{{\tau_i\in[0, \tau)}}\int_{\tau_{i-1}\wedge s}^{\tau_i\wedge s}[-f(r,\wt Y_{\zeta_{i-1}}(r),\wt Z_{\zeta_{i-1}}(r),\zeta_{i-1})dr+\wt Z_{\zeta_{i-1}}(r)dW(r)-d\wt K_{\zeta_{i-1}}(r)]\\
& &+\mathop{\sum_{i\ge1}}_{{\tau_i\in[0, \tau)}}[\wt Y_{\zeta_i}(\tau_i)-\wt  Y_{\zeta_{i-1}}(\tau_i)]\1_{[\tau_i,\tau]}(s)\\
&=&\int_0^s[-f(r,\wt Y^{\alpha_\cd}(r),\wt Z^{\alpha_\cd}(r),\alpha_r)dr+\wt Z^{\alpha_\cd}(r)\,dW(r)- d\wt K^{\alpha(\cdot)}(r)]+\tilde{A}^{\alpha(\cdot)}(s)-A^{\alpha(\cdot)}(s),
\end{eqnarray*}
where
\be\label{UAtilde}\begin{array}{l}
\wt A^{\alpha(\cdot)}(s)=\mathop{\sum_{i\ge1}}_{{\tau_i\in[0, \tau)}}\left[\wt Y_{\zeta_i}(\tau_i)+C_{\zeta_{i-1},\zeta_i}-\wt Y_{\zeta_{i-1}}(\tau_i)\right]\1_{[\tau_i,\tau]}(s),
\end{array}\ee
which is increasing since we have
$$\wt Y(t)\in \bar{Q},\quad \forall t\ge 0.$$
Thus it implies that $(\wt Y^{\alpha_\cd},\wt
Z^{\alpha_\cd})$ is a solution of the following BSDE:
\begin{eqnarray}\label{bsdetilde}
\wt Y^{\alpha_\cd}(s)
&=&g(X^{\alpha(\tau)}(\tau))+A^{\alpha(\cdot)}(\tau)-A^{\alpha(\cdot)}(s)+[(\wt K^{\alpha(\cdot)}(\tau)+\wt A^{\alpha(\cdot)}(\tau))- (\wt K^{\alpha(\cdot)}(s)+\wt A^{\alpha(\cdot)}(s))]\nonumber \\
& &+\int_{s\wedge \tau}^\tau f(r,\wt Y^{\alpha_\cd}(r),\wt Z^{\alpha_\cd}(r),\alpha(r))dr-\int_{s\wedge \tau}^\tau\wt Z^{\alpha_\cd}(r)dW(r), \quad s \ge 0 .
\end{eqnarray}

Since both $\wt K^{\alpha(\cdot)}$ and $\wt A^{\alpha(\cdot)}$ are
increasing càdlàg  processes, from the comparison theorem for multi-dimensional infinite horizon  BSDEs of Shi and Zhang  \cite[Theorem 6]{ZS10}
we conclude that
$$\wt Y^{\alpha(\cdot)}(0)\ge U^{\alpha(\cdot)}(0),$$
which implies that
$$\wt Y_i(0)\ge U^{\alpha(\cdot)}(0).$$

The rest of the proof is similar to the proof of Theorem 3.1 in \cite{HT10}.
\ep

\section{Application to optimal switching problem with unbounded stopping time}
In this section we make the link between the optimal switching problem and the infinite horizon multi-dimensional reflected BSDEs studied previously.
We assume that $C$ satisfies Hypothesis \ref{k'}, and we assume also the following hypothesis.

\begin{Hypothesis}\label{Hl} $ $
\begin{enumerate}[(i)]
\item $l(\cd, 0):= (l(\cd, 0, 1), \cdots,  l(\cd, 0, m))^T$ belongs to $M^2$.
\item  For any $t, x, y,$ and $i\in \Lambda$ there exist $\mu_3 \in \R$ and one positive deterministic bounded function $u(t)$, such that
\be
\langle y, l(t, x, i)\rangle &\le& \mu_3 |y|^2+ u(t)|x|^2, \quad \P-\mbox{a.s.,}
\ee
and $\int_0^\infty u(t)dt<\infty$, $\int_0^\infty u^2(t)dt<\infty$.
\item  For any $t, x$ and $i\in \Lambda$ we have
\b*|l(t,x,i)-l(t,x',i)| &\leq& u(t)|x-x'|.\e*
\item $\sigma$ is invertible and $\sigma^{-1}$ is bounded.
\item $b$ is bounded.
\end{enumerate}
\end{Hypothesis}
Under Hypothesis (H1) and  assumptions $(iv)$-$(v)$,  the following stochastic differential equation:
\b*
dX_t &=&\sigma(t,X)dW_t, \quad X_0=x\in \R^d, \quad \quad t \ge 0,
\e*
has a unique solution.
Identically as in the previous section, a switching strategy $\alpha_\cd$ consists in a sequence $\alpha_\cd:=(\tau_k,\zeta_k)_{k\geq 1}$, where $(\tau_k)_{k\geq 1}$ is an increasing sequence of $\F$-stopping times (i.e $\tau_0=0$, $\tau_k\leq \tau_{k+1}$ and $\lim_{k\rightarrow \infty}\tau_k=\tau$), and $\zeta_{k}$ are $\Fc_{\tau_k}$-measurable random variables valued in $\Lambda$. To a strategy $\alpha_\cd=(\tau_k,\zeta_k)_{k\geq 1}$ and an initial regime $i_0$, we associate
the state process $(\alpha_t)_{t\leq \tau}$ defined by
\b*
\alpha_t&:=&\sum_{k\geq 0} \zeta_{k}\mathbf{1}_{[\tau_k,\tau_k+1)}(t), \quad t \ge 0,
\e*
with $\tau_0=0$ and $\zeta_{0}=i_0$. We denote $\Ac$ the set of admissible strategies and $\Ac_i$ the subset of strategies starting from state $i \in \Lambda$ at time $0$:
\b*
\Ac_i&:=&\left\{\alpha_\cdot \in \Ac: \alpha_0=i \mbox{ and } \E^{\alpha_.}\left[\sum_{k\geq1}C_{\zeta_{k-1},\zeta_k}\right]<\infty \right\},
\e*
where $\E^{\alpha_.}$ denotes the expectation w.r.t the probability $P^{\alpha_.}$, defined for each $\alpha_. \in \Ac_i$ on $(\Omega,\Fc)$ by:
\be\label{ProbabilityUnderAlpha}
\Frac{dP^{\alpha_.}}{dP}=exp\left\{ \int_0^\tau b(s,X^{\alpha_.}(s),\alpha_s)dW_s - \frac12 \int_0^\tau |b(s,X^{\alpha_.}(s),\alpha_s|^2ds \right\}.
\ee
From the assumptions on $\sigma$ and $b$, and according to Girsanov's theorem, the process
\b*
B^{\alpha_.}_t&=&B_t - \int_0^{t} b(s,X^{\alpha_.}(s),\alpha_s)ds, \quad  t\geq 0,
 \e*
 is a Brownian motion on $(\Omega,\Fc,P^{\alpha_.})$. Moreover, for each $\alpha_. \in \Ac_i$, $X^{\alpha_.}$ is a weak solution of:
\be \label{weak-sde}
dX^{\alpha_.}_t=\sigma(t,X^{\alpha_.}(t),\alpha_.)dW^{\alpha_.}_t+b(t,X^{\alpha_.}(t),\alpha_t)dt, \quad X^{\alpha_.}_0=x, \quad \quad t\geq 0.
\ee

Let $(P^{\alpha_.},B^{\alpha_.},X^{\alpha_.})$ be a weak solution of SDE (\ref{weak-sde}), associated with the admissible switching strategy $\alpha_. \in \Ac_i$. We consider the total profit at horizon $\tau$ defined by
\b*
J(\alpha_.)=\E^{\alpha_.}\left[ g(X^{\alpha_.}(\tau))+\int_0^\tau l(s, X^{\alpha_.}(s), \alpha_s)ds + \Sum_{i\ge 1} C_{\alpha_{i-1}, \alpha_i}  \right].
\e*
The switching problem is to maximize the profit $J(\alpha_.)$ over $\alpha_. \in \Ac_i$, subject to the state equation (\ref{weak-sde}), which consists in finding an optimal strategy $\alpha^*_.\in \Ac_i$ such that
\b*
J(\alpha^*_.)&=&\sup_{\alpha_\cd\in \Ac_i}J(\alpha_.).
\e*
We define $f$ as follows: $\forall (t,x,z,i)\in \R^+\times \R\times \R^d \times \Lambda$,
\be\label{NewDefinitionf}
f(t,x,z,i):=l(t,x,i)+\langle z,b(t,x,i)\rangle.
\ee
Under Hypothesis \ref{Hl}, and the expression (\ref{NewDefinitionf}),  the following RBSDE:
\be\label{rbsde-switching}
\left\{\begin{array}{l}
Y_i(t)=g(X_i(\tau))+\Int_{t\wedge \tau}^\tau f(s,X_{i}(s),Z_i(s),i) ds+\int_{t\wedge \tau}^\tau dK_i(s)-\int_{t\wedge \tau}^\tau Z_i(s)\, dW(s),\vspace{1mm}\\
Y_i(t)\ge\displaystyle \Max_{j\in \Ic}\{Y_j(t)-C_{i,j}(t)\}, \vspace{1mm}\\
\Int_0^\tau\left(Y_i(s)-\Max_{j\not=i}\{Y_j(s)-C_{i,j}\}\right)dK_i(s)=0,
\end{array}
\right.\ee
has a unique solution $(Y,Z,K)\in S^2\times M^2\times A^2$, thanks to Theorems \ref{existence} and \ref{uniqueness}.\\

Now we give the main result of this section:

\begin{Theorem}
Let $\alpha^*_.=(\tau_{n}^{*},\zeta_{n}^{*})_{n\geq 0}$ be the strategy given by $(\tau_{0}^{*},\zeta_{0}^{*})=(0,i_0)$ with $i_0\in \Lambda$ and defined recursively, for $n\geq 1$, by

\be\begin{array}{l}\label{tau*}
\tau_{n}^{*}= \inf\left\{s\geq \tau_{n-1}^*; \quad Y_{\zeta_{n-1}^{*}}(s)=\ds\max_{j\in \Lambda^{-\zeta_{n-1}^{*}}}\left(Y_{j}(s)-C_{\zeta_{n-1}^{*},j}\right)\right\} \wedge \tau, \\
$ $\\
\zeta_{n}^{*}\in argmax\, \left\{j; \quad Y_{\zeta_{n-1}^{*}}(s)=\ds\max_{j\in \Lambda^{-\zeta_{n-1}^{*}}}\left(Y_{j}(s)-C_{\zeta_{n-1}^{*},j}\right)\right\},
\end{array}
\ee
where $\Lambda^{-i}:= \Lambda-\{i\}$. \\
Under Hypotheses \ref{k} and \ref{Hl}, the strategy $\alpha^*_.$ is optimal for the switching problem and we have
$$Y_{i_0}(0)=J(\alpha^*_.)=\sup_{\alpha_.\in \Ac_{i_0}}J(\alpha_.).$$
\end{Theorem}
{\bf Proof.}
The proof is performed in two steps.\vspace{2mm}\\
{\bf Step 1.} The strategy $\alpha^*_.$ satisfies $Y_{i_0}(0)=J(\alpha^*_.)$.\\
We consider the reflected BSDE (\ref{rbsde-switching})
$$Y_{i_0}(t)=g(X_{i_0}(\tau))+\Int_{t\wedge \tau}^\tau f(s,X_{i_0}(s),Z_{i_0}(s),{i_0}) ds+\int_{t\wedge \tau}^\tau dK_{i_0}(s)-\int_{t\wedge \tau}^\tau Z_{i_0}(s)\, dW(s).$$
Since $Y_i(0)$ is deterministic, then
\begin{eqnarray*}
Y_{i_0}(0)&=&\E^{\alpha^*_.}\left[g(X_{i_0}(\tau))+\Int_0^\tau f(s,X_{i_0}(s), Z_{i_0}(s),i_0) ds+\int_0^\tau dK_{i_0}(s)-\int_0^\tau Z_{i_0}(s)dW(s)\right] \\
&=&\E^{\alpha^*_.}\left[\Int_{0}^{\tau_{1}^{*}}f(s,X_{i_0}(s),Z_{i_0}(s),i_0) ds+K_{i_0}(\tau_1^*)-\int_{0}^{\tau_{1}^{*}} Z_{i_0}(s)dW(s)+Y_{i_0}(\tau_1^*)\right]\\
&=&\E^{\alpha^*_.}\left[\Int_{0}^{\tau_{1}^{*}}l(s,X_{i_0}(s),i_0) ds+K_{i_0}(\tau_1^*)-\int_{0}^{\tau_{1}^{*}} Z_{i_0}(s)dW^{\alpha^*_.}(s)+Y_{i_0}(\tau_1^*)\right].
\end{eqnarray*}
From the definition of $\tau_1^*$ we know that the process $K_{i_0}(\tau_1^*)$ does not increase between $0$ and $\tau_1^*$ and then $K_{i_0}(\tau_1^*)=0$. On the other hand using the Burkholder-Davis-Gandy's inequality and the assumptions on $b$, we have that $\left(\Int_{0}^{t\wedge \tau} Z_{i_0}(s)dW^{\alpha^*_.}(s), t \ge 0\right)$ is a $P^{\alpha^*_.}$-martingale. Therefore

\b*
Y_{i_0}(0)&=&\E^{\alpha^*_.}\left[\Int_{0}^{\tau_{1}^{*}}l(s,X_{i_0}(s),i_0) ds+Y_{i_0}(\tau_1^*)\right].
\e*
From (\ref{tau*}), we have $Y_{i_0}(\tau_1^*)=Y_{\zeta_1^*}(\tau_1^*)-C_{{i_0},\zeta_1^*}$, therefore
\b*
Y_{i_0}(0)&=&\E^{\alpha^*_.}\left[\Int_{0}^{\tau_{1}^{*}}l(s,X_{i_0}(s),i_0) ds+Y_{\zeta_1^*}(\tau_1^*)-C_{i_0,\zeta_1^*}\right].
\e*
In the same spirit, we repeat this reasoning for $Y_{\zeta_1^*}(\tau_1^*)$. We deduce recursively that
\b*
Y_{i_0}(0)&=&\E^{\alpha^*_.}\left[\sum_{k=1}^{n}\Int_{\tau_{k-1}^{*}}^{\tau_{k}^{*}}l(s,X_{\zeta_k^*}(s),\zeta_k^*) ds+Y_{\zeta_n^*}(\tau_n^*)-\sum_{k=1}^{n}C_{\zeta_{k-1}^*,\zeta_k^*}\right],
\e*
where $\Sum_{k=1}^{n}\Int_{\tau_{k-1}^{*}}^{\tau_k^{*}}l(s,X_{\zeta_k^*}(s),\zeta_k^*)ds=\Int_{0}^{\tau_{n}^{*}}l(s,X^{\alpha_.}(s),\alpha_.)ds$.\\
Then the strategy $\alpha^*_.$ is admissible i.e $E^{\alpha^*_.}[\sum_{k\geq1}C_{\zeta_{k-1}^*,\zeta_k^*}]<+\infty$, because if not, we would have
$Y_{i_0}(0)=-\infty$ which contradicts the assumption $Y_{i_0} \in \Sc^2$.
Thus sending $n$ to infinity, we get that
\b*
Y_{i_0}(0)&=&\E^{\alpha^*_.}\left[\Int_{0}^{\tau}l(s,X^{\alpha_.}(s),\alpha_.)ds-\sum_{k\geq1}C_{\zeta_{k-1}^*,\zeta_k^*}+Y_{\alpha^*_\tau}(\tau)\right],
\e*
where $Y_{\alpha^*_\tau}(\tau)=g(X_\tau)$.  Therefore we obtain that $Y_{i_0}(0)=J(\alpha^*_.)$.
\vspace{2mm}\\
{\bf Step 2.} The strategy $\alpha^*_.$ is optimal.\\
We pick any strategy $\alpha_.=(\tau_{n},\zeta_{n})_{n\geq 0} \in \Ac_{i_0}$, we consider once again the reflected BSDE (\ref{rbsde-switching})
\b*
Y_{i_0}(0)&=&\E^{\alpha_.}\left[g(X_{i_0}(\tau))+\Int_0^\tau f(s,X_{i_0}(s), Z_{i_0}(s),{i_0}) ds+\int_0^\tau dK_{i_0}(s)-\int_0^\tau Z_{i_0}(s)dW(s)\right] \\
&=&\E^{\alpha_.}\left[\Int_{0}^{\tau_{1}}f(s,X_{i_0}(s), Z_{i_0}(s),{i_0}) ds+K_{i_0}(\tau_1)-\int_{0}^{\tau_{1}} Z_{i_0}(s)dW(s)+Y_{i_0}(\tau_1)\right]\\
&=&\E^{\alpha_.}\left[\Int_{0}^{\tau_{1}}l(s,X_{i_0}(s),{i_0}) ds+K_{i_0}(\tau_1)-\int_{0}^{\tau_{1}} Z_{i_0}(s)dW^{\alpha_.}(s)+Y_{i_0}(\tau_1)\right].
\e*
On the one hand, we know that $Y_{i_0}(\tau_1)\geq Y_{\zeta_1}(\tau_1)-C_{{i_0},\zeta_1}$ and $K_{i_0}(\tau_1)\geq 0$. Moreover $\left(\Int_{0}^{t\wedge \tau} Z_{i_0}(s)dW^{\alpha_.}(s), t\ge 0\right)$ is a $P^{\alpha_.}$-martingale, therefore:
\b*
Y_{i_0}(0)&\geq& \E^{\alpha_.}\left[\Int_{0}^{\tau_{1}}l(s,X_{i_0}(s),{i_0}) ds+Y_{\zeta_1}(\tau_1)-C_{{i_0},\zeta_1}\right].
\e*
Next we replace $Y_{\zeta_1}(\tau_1)$ by its value using the same reasoning, and by proceeding exactly as in step 1, an induction argument leads to
\b*
Y_{i_0}(0)&\geq& \E^{\alpha_.}\left[\sum_{k=1}^{n}\Int_{\tau_{k-1}}^{\tau_{k}}l(s,X_{\zeta_k}(s),\zeta_k) ds+Y_{\zeta_n}(\tau_n)-\sum_{k=1}^{n}C_{\zeta_k-1,\zeta_k}\right].
\e*
Sending $n$ to infinity, since the strategy is admissible, we get
\b*
Y_{i_0}(0)&\ge& \E^{\alpha_.}\left[\Int_{0}^{\tau}l(s,X^{\alpha_.}(s),\alpha_.)ds-\sum_{k\geq1}C_{\zeta_{k-1},\zeta_k}+Y^{\alpha_\tau}(\tau)\right],
\e*
with $Y^{\alpha_\tau}(\tau)=g(X^{\alpha_\tau}_\tau)$.\\

Therefore we obtain that $Y_{i_0}(0)\geq J(\alpha_.).$
The arbitrariness of $\alpha_.$ concludes the proof.
\ep


\begin{thebibliography}{99}

\bibitem{AO06}  K. Akdim, and Y. Ouknine, Infinite horizon reflected backward SDEs with jumps and RCLL obstacle,
{\it Stoch. Anal. Appl.}, 24 (2006), no. 6, 1239-1261.

\bibitem{CL08}   R. Carmona and M. Ludkovski, Pricing asset scheduling flexibility using optimal switching,
{\it Appl. Math. Finance }15 (2008), no. 5-6, 405-447.

\bibitem{C98} Z. Chen, Existence and uniqueness for BSDE with stopping time,
{\it Chinese science bulletin}, Volume 43, N 2 (1998), 96-99.

\bibitem{E10}    B. El Asri, Optimal multi-modes switching problem in infinite horizon,
{\it Stoch. Dyn.} 10 (2010), no. 2, 231-261.


\bibitem{E97}    N. El Karoui, C. Kapoudjian, E. Pardoux, S. Peng, and M. C. Quenez, . Reflected solutions of backward SDE's, and related obstacle problems for PDE's.
{\it Ann. Probab}. 25 (1997), no. 2, 702-737.


\bibitem{HLW99} S. Hamadène, J. P. Lepeltier, and Z.  Wu, Infinite horizon reflected backward stochastic differential equations and applications in mixed control and game problems,
{\it Probab. Math. Statist.} 19 (1999), no. 2, Acta Univ. Wratislav. No. 2198, 211-234.

\bibitem{HT10}  Y. Hu and S. Tang, Multi-dimensional BSDE with oblique reflection and optimal switching,
{\it Probability theory and related fields}, Volume 147, Numbers 1-2 (2010), 89-121.

\bibitem{GP96} A. Gegout-Petit and  E. Pardoux, Equations différentielles rétrogrades réfléchies dans un convexe,
{\it Stochastics and Stochastics Reports}, 57 (1996), 111-128.

\bibitem{M00} M. Meyer, Continuous Stochastic Calculus with Applications to Finance,
{\it Chapman and Hall/CRC}, 1 edition (2000).

\bibitem{P99} E. Pardoux, BSDEs, weak convergence and homogenization of semilinear PDEs,
 {\it Nonlinear analysis, differential equations and control} (Montreal, QC, 1998), 503-549, NATO Sci. Ser. C Math. Phys. Sci., 528, Kluwer Acad. Publ., Dordrecht, 1999.

\bibitem{PTW09}   A. Porchet, N. Touzi and X. Warin, Valuation of power plants by utility indifference and numerical computation,
 {\it Math. Methods Oper. Res.} 70 (2009), no. 1, 47-75.

\bibitem{ZS10} Y. Shi and L. Zhang,  Comparison Theorems of Infinite Horizon Forward-Backward Stochastic Differential Equations,
{\it http://arxiv.org/abs/1005.4139}.

\bibitem{W04}   Z. Wu, Forward-Backward Stochastic Differential Equations with Stopping Time,
{\it Acta Mathematica Scientia}, Ser. B, Vol. 24, No. 1 (2004), 91-99.

\bibitem{Y08}   J. Yin, On solutions of a class of infinite horizon FBSDEs,
{\it Statistics and  Probability Letters}, Volume 78, Issue 15, (2008),2412-2419.


\end{thebibliography}
\end{document}